\documentclass[a4paper]{llncs}
\usepackage[english]{babel}
\usepackage{german}

\usepackage{cite}
\usepackage{amsmath}
\usepackage{amssymb}
\usepackage{enumerate,enumitem}
\usepackage{graphicx,subfigure}

\DeclareGraphicsExtensions{.pdf}
\graphicspath{{figs/}}

\spnewtheorem*{openproblem}{Open Problem}{\bf}{\it}

\renewcommand{\L}{{\normalfont{\textsf{L}}}}

\newcommand{\G}{\rotatebox[origin=c]{180}{\reflectbox{\L}}}
\newcommand{\lE}{{\reflectbox{\L}}}
\newcommand{\eeG}{\rotatebox[origin=c]{180}{{\L}}}

\newcommand{\SEG}{\ensuremath{\operatorname{SEG}}}
\newcommand{\PLANAR}{\ensuremath{\operatorname{PLANAR}}}
\newcommand{\COPLANAR}{\ensuremath{\operatorname{CO-PLANAR}}}
\newcommand{\VPG}{\ensuremath{\operatorname{VPG}}}
\newcommand{\EPG}{\ensuremath{\operatorname{EPG}}}
\newcommand{\STRING}{\ensuremath{\operatorname{STRING}}}
\newcommand{\COCO}{\ensuremath{\operatorname{COCOMP}}}

\begin{document}

\title{Intersection Graphs of L-Shapes and Segments in the Plane\thanks{This work was partially supported by (i)~the DFG ESF EuroGIGA projects COMPOSE and GraDR, 
(ii)~the EPSRC Grant~EP/K022660/1, (iii)~the ANR Project EGOS: ANR-12-JS02-002-01, and (iv)~the PEPS grant EROS. A preliminary version appeared in the Proceedings of MFCS 2014~\cite{US-14}.}}

\author{Stefan Felsner\inst{1} \and Kolja Knauer\inst{2} \and George B. Mertzios\inst{3} \and Torsten Ueckerdt\inst{4}}

\institute{
Institut f\"ur Mathematik, Technische Universit\"at Berlin, Germany.\\
Email: \texttt{felsner@math.tu-berlin.de} 
\and
LIF UMR 7279, Universit\'e Aix-Marseille, CNRS, France.\\
Email: \texttt{kolja.knauer@lif.univ-mrs.fr} 
\and
School of Engineering and Computing Sciences, Durham University, UK.\\
Email: \texttt{george.mertzios@durham.ac.uk} 
\and
Department of Mathematics, Karlsruhe Institute of Technology, Germany.\\ 
Email: \texttt{torsten.ueckerdt@kit.edu}\vspace{-0.2cm}}
\date{\vspace{-5cm}}
\maketitle


\begin{abstract}
 An L-shape is the union of a horizontal and a vertical segment with a common endpoint. These come in four rotations: \L,~\G,~\lE{}~and~\eeG. 
 A $k$-bend path is a simple path in the plane, whose direction changes $k$ times from horizontal to vertical.
 If a graph admits an intersection representation in which every vertex is represented by an \L, an \L{} or \G, a $k$-bend path, or a segment, then this graph is called an $\{\L\}$-graph, $\{\L,\G\}$-graph, $B_k$-\VPG-graph or \SEG-graph, respectively. Motivated by a theorem of Middendorf and Pfeiffer [Discrete Mathematics, 108(1):365--372, 1992], stating that every $\{\L,\G\}$-graph is a \SEG-graph, we investigate several known subclasses of \SEG-graphs and show that they are $\{\L\}$-graphs, or $B_k$-\VPG{}-graphs for some small constant $k$.
 We show that all planar $3$-trees, all line graphs of planar graphs, and all full subdivisions of planar graphs are $\{\L\}$-graphs. Furthermore we show that complements of planar graphs are $B_{17}$-\VPG-graphs and complements of full subdivisions are $B_2$-\VPG-graphs. Here a full subdivision is a graph in which each edge is subdivided at least once.\newline

\end{abstract}

\section{Introduction and Motivation}
 \label{sec:introduction}
 
 A \textbf{segment intersection graph}, \SEG-graph for short, is a graph that can be represented as follows. Vertices correspond to straight-line segments in the plane and two vertices are adjacent if and only if the corresponding segments intersect. Such representations are called \emph{\SEG-representations} and, for convenience, the class of all \SEG-graphs is denoted by \SEG. \SEG-graphs are an important subject of study strongly motivated from an algorithmic point of view. Indeed, having an intersection representation of a graph (in applications graphs often come along with such a given representation) may allow for designing better or faster algorithms for optimization problems that are hard for general graphs, such as finding a maximum clique in interval graphs.
 
 More than 20 years ago, Middendorf and Pfeiffer~\cite{MP92}, considered intersection graphs of \textbf{axis-aligned L-shapes} in the plane, where an axis-aligned L-shape is the union of a horizontal and a vertical segment whose intersection is an endpoint of both. In particular, L-shapes come in four possible rotations: \L, \G, \lE, and \eeG. For a subset $X$ of these four rotations, e.g., $X = \{\L\}$ or $X = \{\L,\G\}$, we call a graph an \emph{$X$-graph} if it admits an \emph{$X$-representation}, i.e., vertices can be represented by L-shapes from $X$ in the plane, each with a rotation from $X$, such that two vertices are adjacent if and only if the corresponding L-shapes intersect. Similarly to \SEG{}, we denote the class of all $X$-graphs by $X$. The question if an intersection representation with polygonal paths or pseudo-segments can be \emph{stretched} into a \SEG-representation is a classical topic in combinatorial geometry and Oriented Matroid Theory. Middendorf and Pfeiffer prove the following 
interesting relation between intersection graphs of segments and L-shapes.
 
 \begin{theorem}[Middendorf and Pfeiffer~\cite{MP92}]\label{thm:stretch-L-G}
  Every $\{\L,\G\}$-representation has a combinatorially equivalent \SEG-representation.
 \end{theorem}
 
 This theorem is best-possible in the sense that there are examples of $\{\L, \eeG\}$-graphs which are no \SEG-graphs~\cite{CJKV12,MP92}, i.e., such $\{\L, \eeG\}$-representations cannot be stretched.
 We feel that Theorem~\ref{thm:stretch-L-G}, which of course implies that $\{\L,\G\} \subseteq \SEG$, did not receive a lot of attention in the active field of \SEG-graphs. In particular, one could use Theorem~\ref{thm:stretch-L-G} to prove that a certain graph class $\mathcal{G}$ is contained in \SEG{} by showing that $\mathcal{G}$ is contained in $\{\L,\G\}$. For example, very recently Pawlik \emph{et al.}~\cite{Paw-13} discovered a class of triangle-free \SEG-graphs with arbitrarily high chromatic number, disproving a famous conjecture of Erd\H{o}s~\cite{Gya-87}, and it is in fact easier to see that these graphs are $\{\L\}$-graphs than to see that they are \SEG-graphs. To the best of our knowledge, the stronger result $\mathcal{G} \subseteq \{\L,\G\}$ has never been shown for any non-trivial graph class $\mathcal{G}$. In this paper we initiate this research direction. We consider several graph classes which are known to be contained in \SEG{} and show that they are actually contained in $\{\L\}$, which 
is a proper subclass of $\{\L,\G\}$~\cite{CJKV12}.
  
Whenever a graph is not known (or known not) to be an intersection graph of segments or axis-aligned L-shapes, one often considers natural generalizations of these intersection representations. Asinowski \textit{et al.}~\cite{ACGLLS12} introduced \textbf{intersection graphs of axis-aligned $k$-bend paths} in the plane, called $B_k$-\VPG{}-graphs. 
An (axis-aligned) $k$-bend path is a simple path in the plane, whose direction changes $k$ times from horizontal to vertical.
Clearly, $B_1$-\VPG{}-graphs are precisely intersection graphs of all four L-shapes; the union of $B_k$-\VPG-graphs for all $k \geq 0$ is exactly the class \STRING{} of intersection graphs of simple curves in the plane~\cite{ACGLLS12}. Now if a graph $G \notin \SEG$ is a $B_k$-\VPG{}-graph for some small $k$, then one might say that $G$ is ``not far from being a \SEG-graph''.

\bigskip

\noindent
\textbf{Our Results and Related Work.}{\ \\}
%
\noindent
 Let us denote the class of all planar graphs by \PLANAR. A recent celebrated result of Chalopin and Gon\c{c}alves~\cite{Cha-09} states that $\PLANAR \subset \SEG$, 
which was conjectured by Scheinerman~\cite{Sch-84} in 1984. However, their proof is rather involved and there is not much control over the kind of \SEG-representations. 
Here we give an easy proof for a non-trivial subclass of planar graphs, namely \emph{planar 3-trees}. 
A \emph{3-tree} is an edge-maximal graph of treewidth 3. Every 3-tree can be built up starting from the clique $K_{4}$ and adding new vertices, one at a time, whose neighborhood in the so-far 
constructed graph is a triangle. 
 
 
 \begin{theorem}\label{thm:planar-3-trees}
  Every planar $3$-tree is an $\{\L\}$-graph.
 \end{theorem}
 
 
 It remains open to generalize Theorem~\ref{thm:planar-3-trees} to planar graphs of treewidth~$3$ (i.e., subgraphs of planar 3-trees). On the other hand it is easy to see that graphs of treewidth at most $2$ are $\{\L\}$-graphs~\cite{CKU-13}. Chaplick and the last author show in~\cite{CU13} that planar graphs are $B_2$-\VPG{}-graphs, improving on an earlier result of Asinowski \textit{et al.}~\cite{ACGLLS12}. In~\cite{CU13} it is also conjectured that $\PLANAR \subset \{\L\}$, which with Theorem~\ref{thm:stretch-L-G} would imply the main result of~\cite{Cha-09}, i.e., $\PLANAR \subset \SEG$.
 
 
 Considering line graphs of planar graphs, one easily sees that these graphs are \SEG-graphs. Indeed, a straight-line drawing of a planar graph $G$ can be interpreted as a \SEG-representation of the line graph $L(G)$ of $G$, which has the edges of $G$ as its vertices and pairs of incident edges as its edges. We prove the following strengthening result.
 
 \begin{theorem}\label{thm:planar-linegraphs}
  The line graph of every planar graph is an $\{\L\}$-graph.
 \end{theorem}
 

 Kratochv{\'{\i}}l and Kub{\v{e}}na~\cite{Kra-98} consider the class of complements of planar graphs (co-planar graphs), \COPLANAR{} for short. They show that every graph in \COPLANAR{} is an intersection graph of convex sets in the plane, and ask whether $\COPLANAR{} \subset \SEG$. As the \textsc{Independent Set Problem} in planar graphs is known to be NP-complete~\cite{GJ-79}, \textsc{Max Clique} is NP-complete for any graph class $\mathcal{G} \supseteq \COPLANAR{}$, e.g., intersection graphs of convex sets. Indeed, the longstanding open question whether \textsc{Max Clique} is NP-complete for \SEG{}~\cite{Kra-94} has recently been answered affirmatively by Cabello, Cardinal and Langerman~\cite{Cab-13} by showing that every planar graph has an even subdivision whose complement is a \SEG-graph. The subdivision is essential in the proof of~\cite{Cab-13}, as it still remains an open problem whether $\COPLANAR \subset \SEG$~\cite{Kra-98}. The largest subclass of \COPLANAR{} known to be in \SEG{} is the 
class of complements 
of partial $2$-trees~\cite{Fra-12}. Here we show that all co-planar graphs are ``not far from being \SEG-graphs''.
 
 \begin{theorem}\label{thm:co-planar}
  Every co-planar graph is a $B_{17}$-VPG-graph.
 \end{theorem}
 
 Theorem~\ref{thm:co-planar} implies that \textsc{Max Clique} is NP-complete for $B_k$-\VPG-graphs with $k \geq 17$. On the other hand, the \textsc{Max Clique} problem for $B_0$-\VPG-graphs can be solved in polynomial time, while \textsc{Vertex Colorability} remains NP-complete but allows for a $2$-approxima\-tion~\cite{ACGLLS12}. Middendorf and Pfeiffer~\cite{MP92} show that the complement of any \emph{even subdivision} of any graph, i.e., every edge is subdivided with a non-zero even number of vertices, is an $\{\L,\eeG\}$-graph. This implies that \textsc{Max Clique} is NP-complete even for $\{\L,\eeG\}$-graphs.
 
 \smallskip

 We consider \emph{full subdivisions} of graphs, that is, a subdivision $H$ of a graph $G$ where each edge of $G$ is subdivided at least once. It is not hard to see that a full subdivision $H$ of $G$ is in \STRING{} if and only if $G$ is planar, and that if $G$ is planar, then $H$ is actually a \SEG-graph. Here we show that this can be further strengthened, namely that $H$ is in an $\{\L\}$-graph. 
Moreover, we consider the complement of a full subdivision $H$ of an arbitrary graph $G$, which is in \STRING{} but not necessarily in \SEG. Here, similar to the result of Middendorf and Pfeiffer~\cite{MP92} on even subdivisions we show that such a graph $H$ is ``not far from being \SEG-graph''.
 
 \begin{theorem}\label{thm:co-subdivision}
  Let $H$ be a full subdivision of a graph $G$.
  \begin{enumerate}[label = (\roman*)]
   \item If $G$ is planar, then $H$ is an $\{\L\}$-graph.
   \item If $G$ is any graph, then the complement of $H$ is a $B_2$-\VPG-graph.
  \end{enumerate}
 \end{theorem}
 
\begin{figure}[htb]
  \centering
  \includegraphics[scale=0.9]{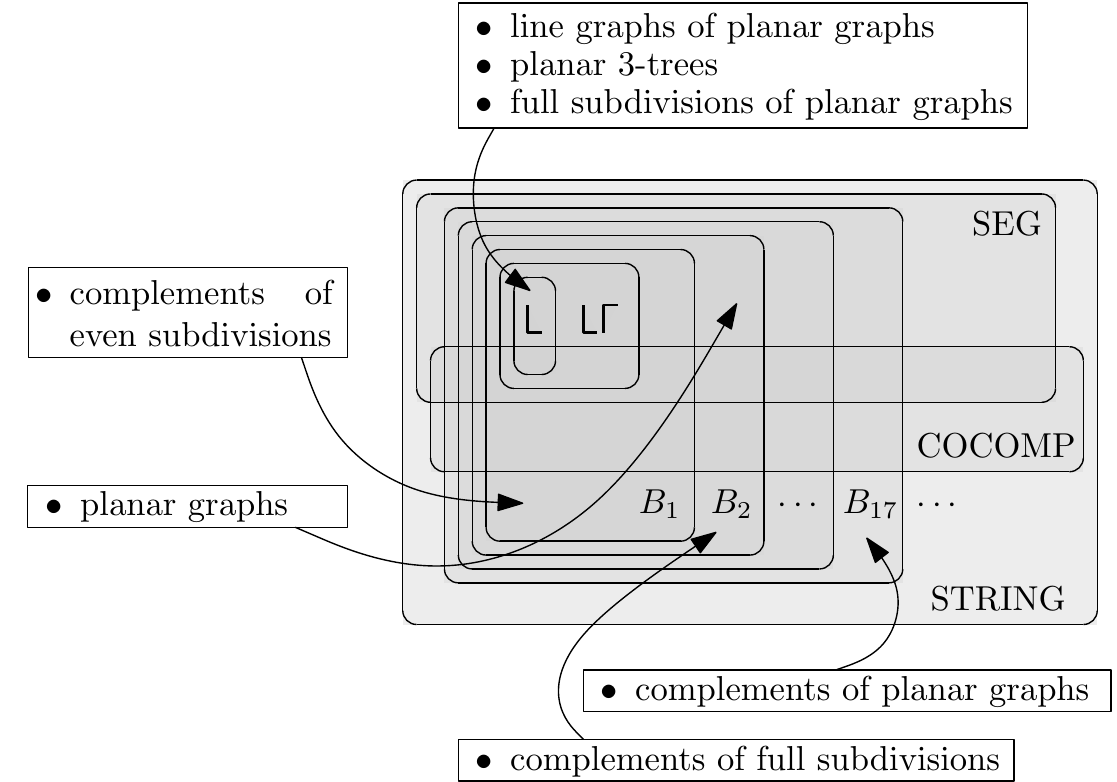}
  \caption{Graph classes considered in this paper.}
  \label{fig:graph-classes}
 \end{figure}

The graph classes considered in this paper are illustrated in Figure~\ref{fig:graph-classes}. 
We shall prove Theorems~\ref{thm:planar-3-trees},~\ref{thm:planar-linegraphs},~\ref{thm:co-planar} and~\ref{thm:co-subdivision} in Sections~\ref{sec:planar-3-trees},~\ref{sec:planar-linegraphs},~\ref{sec:co-planar} and~\ref{sec:co-subdivision}, respectively, and conclude with some open questions in Section~\ref{sec:conclusions}. 


\bigskip

\noindent
\textbf{Related Representations.}{\ \\}
 \noindent
 In the context of \emph{contact representations}, where distinct segments or $k$-bend paths may not share interior points, it is known that every contact \SEG-representation has a combinatorially equivalent contact $B_1$-\VPG-re\-pre\-sen\-ta\-tion, but not vice versa~\cite{KUV-13}. Contact \SEG-graphs are exactly planar Laman graphs and their subgraphs~\cite{Fra-07}, which includes for example all triangle-free planar graphs. Very recently, contact $\{\L\}$-graphs have been characterized~\cite{CKU-13}. Necessary and sufficient conditions for stretchability of a contact system of pseudo-segments are known~\cite{deFraysseix2007stretching,aerts2012triangle}.
 
 \smallskip
 
%
%

 
 Let us also mention the closely related concept of \emph{edge}-intersection graphs of paths in a grid (\EPG-graphs) introduced by Golumbic \textit{et al.}~\cite{Gol-09}. There are some notable differences, starting from the fact that \emph{every} graph is an \EPG-graph~\cite{Gol-09}. Nevertheless, analogous questions to the ones posed about \VPG-representations of \STRING-graphs are posed about \EPG-representations of general graphs. In particular, there is a strong interest in finding representations using paths with few bends, see~\cite{Hel-12} for a recent account.

\section{Proof of Theorem~\ref{thm:planar-3-trees}}
\label{sec:planar-3-trees}

\begin{proof}


Let $G$ be a plane $3$-tree with a fixed plane embedding. 
  We construct an $\{\L\}$-representation of $G$ satisfying the additional property that for every inner triangular face $\{a,b,c\}$ of $G$ there exists a subset of the plane, called the \emph{private region of the face}, that intersects only the \L-paths for $a$, $b$ and $c$, and no other \L-path.
More precisely, a private region of $\{a,b,c\}$ is an axis-aligned rectilinear polygon having one of the shapes depicted in Figure~\ref{fig:private_regions}, such that the \L-paths for $a$, $b$ and $c$ intersect the polygon as shown in figure.
  
  \begin{figure}[htb]
   \centering
   \subfigure[\label{fig:private_regions}]{
    \includegraphics{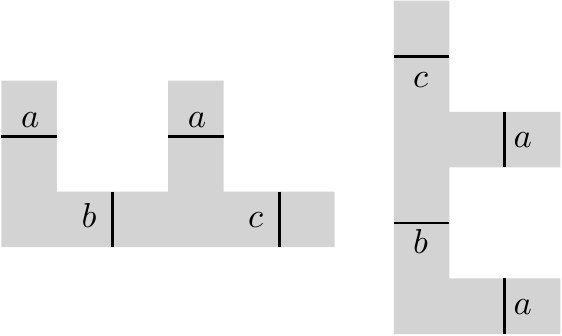}
   }
   \hspace{4em}
   \subfigure[\label{fig:induction_base}]{
    \includegraphics{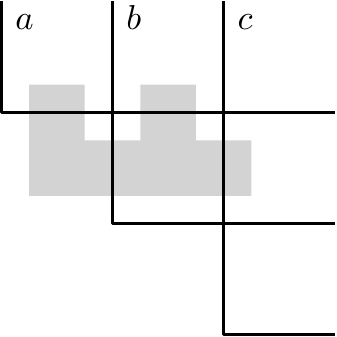}
   }
   \caption{\subref{fig:private_regions} The two possible shapes of a private region for inner facial triangle $\{a,b,c\}$. \subref{fig:induction_base} An $\{\L\}$-representation of the plane $3$-tree on three vertices together with a private region for the only inner face.}
  \end{figure}
  
  Indeed, we prove the following stronger statement by induction on the number of vertices in $G$.
  
  \begin{claim}
   Every plane $3$-tree admits an $\{\L\}$-representation together with a private region for every inner face, such that the private regions for distinct faces are disjoint.
  \end{claim}

  As induction base ($|V(G)|=3$) consider the graph $G$ consisting only of the triangle $\{a,b,c\}$. Then there is an essentially unique $\{\L\}$-representation of $G$ and it is not difficult to find a private region for the unique inner face of $G$. We refer to Figure~\ref{fig:induction_base} for an illustration.
  
  Now let us assume that $|V(G)|\geq 4$. Because $G$ is a $3$-tree there exists an inner vertex $v$ of degree exactly three. In particular, the three neighbors $a,b,c$ of $v$ form an inner facial triangle in the plane $3$-tree $G' = G \setminus v$. By induction $G'$ admits an $\{\L\}$-representation with a private region for each inner face so that distinct private regions are disjoint. 
  
  Consider the private region $R$ for $\{a,b,c\}$. By flipping the plane along the main diagonal if necessary, we can assume without loss of generality that $R$ has the shape shown in the left of Figure~\ref{fig:private_regions}. (Note that such a flip does not change the type of the \L-paths.) Now we introduce an \L-path for vertex $v$ completely inside $R$ as depicted in Figure~\ref{fig:induction_step}. Since $R$ does not intersect any other \L-path this is an $\{\L\}$-representation of $G$.
  
  \begin{figure}[htb]
   \centering
   \subfigure[\label{fig:induction_step}]{
    \includegraphics{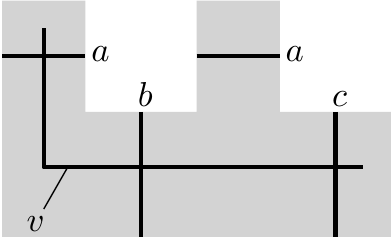}
   }
   \hspace{3em}
   \subfigure[\label{fig:new_regions}]{
    \includegraphics{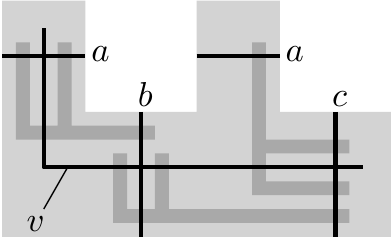}
   }
   \caption{\subref{fig:induction_step} Introducing an \L-shape for vertex $v$ into the private region for the triangle $\{a,b,c\}$. \subref{fig:new_regions} Identifying a pairwise disjoint private regions for the facial triangles $\{a,b,v\}$, $\{a,c,v\}$ and $\{b,c,v\}$.}
  \end{figure}
 
  Finally we identify three private regions for the three newly created inner faces $\{a,b,v\}$, $\{a,c,v\}$ and $\{b,c,v\}$. This is shown in Figure~\ref{fig:new_regions}. Since these regions are pairwise disjoint and completely contained in the private region for $\{a,b,c\}$ we have identified a private region for every inner face so that distinct regions are disjoint. (Note that $\{a,b,c\}$ is not a facial triangle in $G$ and hence does not need a private region.) This proves the claim and thus concludes the proof of the theorem.\qed
\end{proof}

\section{Proof of Theorem~\ref{thm:planar-linegraphs}}
\label{sec:planar-linegraphs}
 
 \begin{proof}
  Without loss of generality let $G$ be a maximally planar graph with a fixed plane embedding. (Line graphs of subgraphs of $G$ are induced subgraphs of $L(G)$.) Then $G$ admits a so-called \emph{canonical ordering} --first defined in~\cite{deF-90}--, namely an ordering $v_1,\ldots,v_n$ of the vertices of $G$ such that
  \begin{itemize}
   \item Vertices $v_1,v_2,v_n$ form the outer triangle of $G$ in clockwise order. (We draw $G$ such that $v_1,v_2$ are the highest vertices.) 
   \item For $i=3,\ldots,n$ vertex $v_i$ lies in the outer face of the induced embedded subgraph $G_{i-1} = G[v_1,\ldots,v_{i-1}]$. Moreover, the neighbors of $v_i$ in $G_{i-1}$ form a path on the outer face of $G_{i-1}$ with at least two vertices.
  \end{itemize}
  We shall construct an $\{\L\}$-representation of $L(G)$ along a fixed canonical ordering $v_1,\ldots,v_n$ of $G$. For every $i=2,\ldots,n$ we shall construct an $\{\L\}$-representation of $L(G_i)$ with the following additional properties.
  
  For every outer vertex $v$ of $G_i$ we maintain an auxiliary bottomless rectangle $R(v)$, i.e., an axis-aligned rectangle with bottom-edge at $-\infty$, such that: 
  \begin{itemize} 
  \item $R(v)$ intersects the horizontal segments of precisely those rectilinear paths for edges in $G_i$ incident to $v$. 
  \item $R(v)$ does not contain any bends or endpoints of any path for an edge in $G_i$ and does not intersect any $R(w)$ for $w \neq v$.
  \item the left-to-right order of the bottomless rectangles matches the order of vertices on the counterclockwise outer $(v_1,v_2)$-path of $G_i$.
  \end{itemize}
  The bottomless rectangles act as placeholders for the upcoming vertices of $L(G)$. Indeed, all upcoming intersections of paths will be realized inside the corresponding bottomless rectangles.
  For $i=2$, the graph $G_i$ consist only of the edge $\{v_1,v_2\}$. Hence an $\{\L\}$-representation of the one-vertex graph $L(G_2)$ consists of only one L-shape and two disjoint bottomless rectangles $R(v_1)$, $R(v_2)$ intersecting its horizontal segment.
  
  \medskip
  
  For $i \geq 3$, we start with an $\{\L\}$-representation of $L(G_{i-1})$. Let $(w_1,\ldots,w_k)$ be the counterclockwise outer path of $G_{i-1}$ that corresponds to the neighbors of $v_i$ in $G_{i-1}$. The corresponding bottomless rectangles $R(w_1),\ldots,R(w_k)$ appear in this left-to-right order. See Figure~\ref{fig:planar-line} for an illustration. For every edge $\{v_i, w_j\}$, $j = 1,\ldots,k$ we define an L-shape $P(v_iw_j)$ whose vertical segment is contained in the interior of $R(w_j)$ and whose horizontal segment ends in the interior of $R(w_k)$. Moreover, the upper end and lower end of the vertical segment of $P(v_iw_j)$ lies on the top side of $R(w_j)$ and below all \L-shapes for edges in $G_{i-1}$, respectively. Finally, the bend and right end of $P(v_iw_j)$ is placed above the bend of $P(v_iw_{j+1})$ and to the right of the right end of $P(v_iw_{j+1})$ for $j = 1,\ldots,k-1$, see Figure~\ref{fig:planar-line}.
  
  
  \begin{figure}[htb]
   \centering
   \includegraphics[scale=0.85]{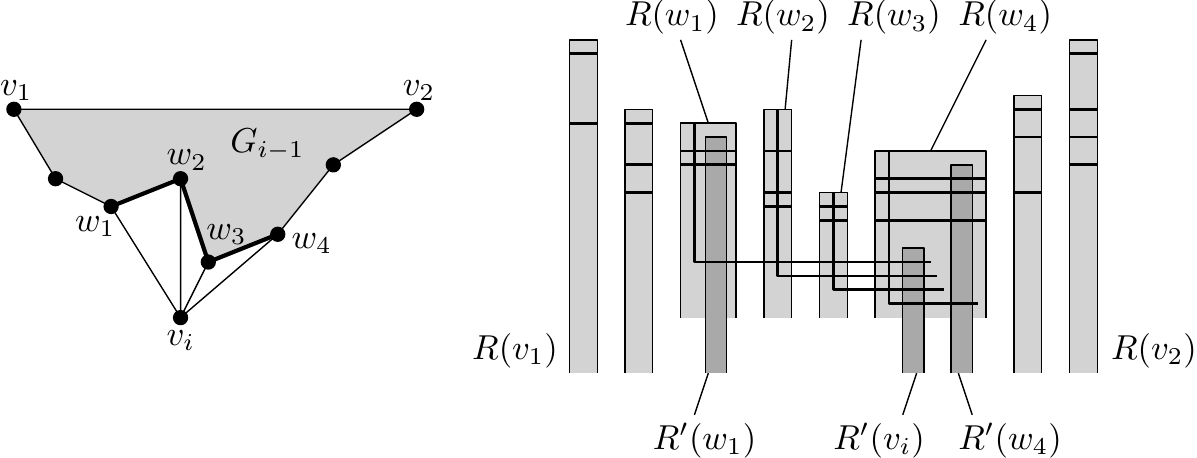}
   \caption{Along a canonical ordering a vertex $v_i$ is added to $G_{i-1}$. For each edge between $v_i$ and a vertex in $G_{i-1}$ an \L-shape is introduced with its vertical segment in the corresponding bottomless rectangle. The three new bottomless rectangles $R'(w_1),R'(v_i),R'(w_k)$ are highlighted.}
   \label{fig:planar-line}
  \end{figure}
  

  It is straightforward to check that this way we obtain an $\{\L\}$-representation of $L(G_i)$. So it remains to find a set of bottomless rectangles, one for each outer vertex of $G_i$, satisfying our additional property. We set $R'(v) = R(v)$ for every $v \in V(G_i) \setminus \{v_i,w_1,\ldots,w_k\}$ since these are kept unchanged. Since $R(w_1)$ and $R(w_k)$ are not valid anymore, we define a new bottomless rectangle $R'(w_1) \subset R(w_1)$ such that $R'(w_1)$ is crossed by all horizontal segments that cross $R(w_1)$ and additionally the horizontal segment of $P(v_iw_1)$. Similarly, we define $R'(w_k) \subset R(w_k)$. And finally, we define a new bottomless rectangle $R'(v_i) \subset R(w_k)$ in such a way that it is crossed by the horizontal segments of exactly $P(v_iw_1),\ldots,P(v_iw_k)$. Note that for $1 < j < k$ the outer vertex $w_j$ of $G_{i-1}$ is not an outer vertex of $G_i$. Then $\{R'(v) \mid v \in v(G_i)\}$ has the desired property. See again Figure~\ref{fig:planar-line}.\qed
 \end{proof}

\section{Proof of Theorem~\ref{thm:co-planar}}
\label{sec:co-planar}
 
\begin{figure}[htb]
   \centering
   \includegraphics[width=\textwidth]{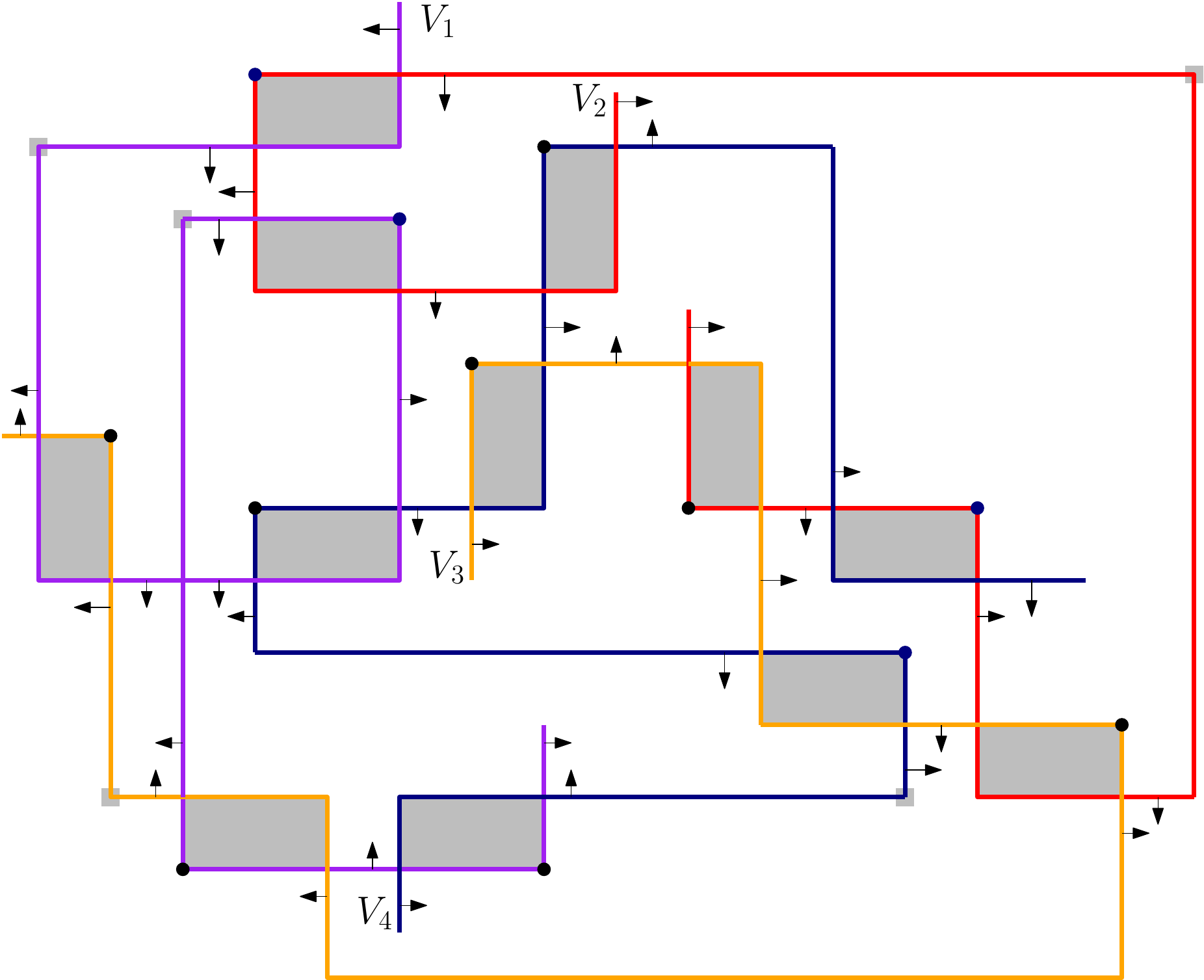}
   \caption{Interconnecting the \VPG-representations by adding at most three bends for each vertex. The set of paths corresponding to color class $V_i$ is indicated by a single path labeled $V_i$, $i=1,2,3,4$. Little gray boxes on bends indicate that the bend is used to reverse the order of coordinates within the class.}
   \label{fig:co-planar-connect}
  \end{figure}

 \begin{proof}
  Let $G = (V,E)$ be any planar graph. We shall construct a $B_k$-VPG-representation of the complement $\bar{G}$ of $G$ for some constant $k$ that is independent of $G$. Indeed, $k=17$ is enough. To find the \VPG-representation we make use of two crucial properties of $G$: A) $G$ is $4$-colorable and B) $G$ is $5$-degenerate. Indeed, our construction gives a $B_{2d+7}$-VPG-representation for the complement of any $4$-colorable $d$-degenerate graph. Here a graph is called \emph{$d$-degenerate} if it admits a vertex ordering such that every vertex has at most $d$ neighbors with smaller index.
  
  Consider any $4$-coloring of $G$ with color classes $V_1,V_2,V_3,V_4$. Further let $v_1,\ldots,v_n$ be an order of the vertices of $V$ witnessing the degeneracy of $G$, i.e., for each $v_i$ there are at most $5$ neighbors $v_j$ of $v_i$ with $j < i$. We call these neighbors the \emph{back neighbors of $v_i$}. 
%
  Further consider the axis-aligned rectangle $R = [0,2(n+1)]\times[0,n+1]$. 
%
%
%
%
%

  Consider any ordered pair of color classes, say $(V_1, V_2)$ and the axis-aligned rectangle $R = [0,2(n+1)]\times[0,n+1]$. We will represent all non-edges of the form $(v_i,v_j)$ with $i<j$, $v_i\in V_1$, and $v_j\in V_2$ in $R$. 
  We define a monotone increasing path $Q(v)$ for each $v \in V_1\cup V_2$ as follows. See Figure~\ref{fig:co-planar} for an illustration.
  
  
  \begin{figure}[htb]
   \centering
   \includegraphics[width=\textwidth]{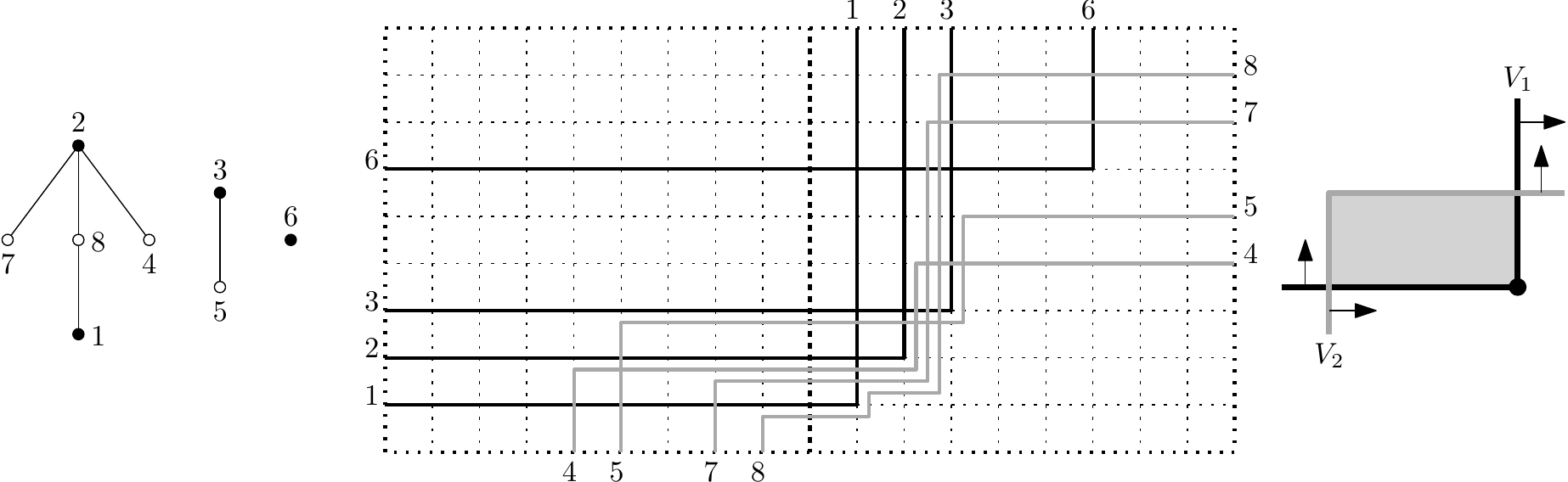}
    \caption{Left: The induced subgraph $G[V_1, V_2]$ for a pair of color classes ($V_1$ (black), $V_2$ (white)) of a planar graph $G$. Middle: A \VPG-representation of all non-edges  of the form $(v_i,v_j)$ with $i<j$, $v_i\in V_1$, and $v_j\in V_2$ in $R$. Right: A schematic representation of the arrangement of paths. T arrows indicated the order of coordinates within $V_1$ and $V_2$. The disk on the bend of the black path indicates, that backward non-edges to $V_1$ are realized.
}
   \label{fig:co-planar}
  \end{figure}
  
  
  \begin{itemize}
   \item For $v_j \in V_1$ let $Q(v_j)$ start on $(0,j)$, go horizontally to $(n+1+j,j)$ and then vertically to $(n+1+j,n+1+j)$. (These are the black paths in the middle of Figure~\ref{fig:co-planar}.)
   
   \item If $v_j \in V_2$ with back neighbors $v_{i_1}, \ldots v_{i_{k'}}\in V_1$ let $Q(v_j)$ start vertically on $(j,0)$, bend to the right at $(j,i_1-\varepsilon)$, bend upwards at $(n+1+i_1+\varepsilon,i_1-\varepsilon)$, and in this way avoid all paths representing back neighbors of $v_j$ in $V_1$. More precisely, $Q(v_j)$ now uses horizontal segments of the form $[(n+1+i_{\ell}+\varepsilon,i_{\ell+1}-\varepsilon),(n+1+i_{\ell+1}+\varepsilon,i_{\ell+1}-\varepsilon)]$ and vertical segments of the form $[(n+1+i_{\ell}+\varepsilon,i_{\ell}-\varepsilon),(n+1+i_{\ell}+\varepsilon,i_{\ell+1}-\varepsilon)]$, for $\ell+1\leq k'$. The last vertical segment of $Q(v_j)$ starts at $[(n+1+i_{k'}+\varepsilon,i_{k'}-\varepsilon)$ and goes to $(n+1+i_{\ell}+\varepsilon,j)$. From there $Q(v_j)$ takes the last horizontal segment to  $(2(n+1),j)$. (These are the gray paths in the middle of Figure~\ref{fig:co-planar}.)

%
%
%
%
%
  \end{itemize}
  
  First, note that in order to avoid collinearities for each $j$ there should be a different $\varepsilon(j)$ in the above construction, which yields a representation on an integer grid of size roughly $2n^2\times n^2$.
  
  Second, note that all paths representing vertices from $V_1$ mutually cross using one bend each, i.e., non-edges between vertices of $V_1$ are represented. 
  Now observe that $\{Q(v) \mid v \in V_1\cup V_2\}$ is a \VPG-representation of all non-edges  of the form $(v_i,v_j)$ with $i<j$, $v_i\in V_1$, and $v_j\in V_2$. Every path $Q(v_j)$ with $v_j\in V_2$ crosses all paths corresponding to $v_i\in V_1$ with $i<j$ except those that are back neighbors of $v_j$. These are avoided at a cost of two bends per back neighbor. A last bend is spent when turning to the right the last time. 
  
  Furthermore, note that the order on the $y$-coordinates of the starting points of the $V_1$-paths and the order on the $x$-coordinates of the last points of the $V_1$-paths are according to their order of indices. Also observe that the order of the $x$-coordinates of the starting points of the $V_2$-paths does not matter for the construction, as long as all of them are in $[0,n+1]$. Last, note that the order on the $y$-coordinates of the last points of the $V_2$-paths is according to their order of indices. We synthesize this situation into a diagram on the right of Figure~\ref{fig:co-planar}, where paths represent the family of paths from one color class, the arrow indicates the order in which the paths within one family are arranged. The disk on the bend of the black path means, that the non-edges from the gray class backwards to the black class are realized.
%
%
  
%
  
  Now we have defined for each ordered pair of color classes $(V_i, V_j)$ a \VPG-re\-pre\-sen\-ta\-tion of the non-edges backward from $V_j$ to $V_i$.
  For every vertex ${v \in V}$ we have defined six $Q$-paths, two for each colors class that $v$ is not in. In total the six $Q$-paths for the same vertex $v$ have at most $6 + 2k \leq 16$ bends, where $k \leq 5$ is the back degree of $v$. It remains to place (translate and rotate) the six representations into non-overlapping positions and to ``connect'' the three $Q$-paths for each vertex in such a way that connections for vertices of different color do not intersect. This can be done with at most three extra bends per path, see Figure~\ref{fig:co-planar-connect}. Finally, note that the first and last segment of every path in the representation can be omitted, yielding the claimed bound.\qed
 \end{proof}
  
  
\section{Proof of Theorem~\ref{thm:co-subdivision}}
\label{sec:co-subdivision}
 
 \begin{proof}
  Let $G$ be any graph and $H$ arise from $G$ by subdividing each edge at least once. Without loss of generality we may assume that every edge of $G$ is subdivided exactly once or twice. Indeed, if an edge $e$ of $G$ is subdivided three times or more, then $H$ can be seen as a full subdivision of the graph $G'$ that arises from $G$ by subdividing $e$ once.
  
  \begin{enumerate}[leftmargin = 1.5em, label=(\roman*)]
   \item Assuming that $G$ is planar, we shall find an $\{\L\}$-representation of $H$ as follows. Without loss of generality $G$ is maximally planar. We consider a bar visibility representation of $G$, i.e., vertices of $G$ are disjoint horizontal segments in the plane and edges are disjoint vertical segments in the plane whose endpoints are contained in the two corresponding vertex segments and which are disjoint from all other vertex segments. Such a representation for a planar triangulation exists e.g. by~\cite{Luc-87}. See Figure~\ref{fig:bar-visibility} for an illustration.
   
   \begin{figure}[htb]
    \centering
    \includegraphics[scale=0.85]{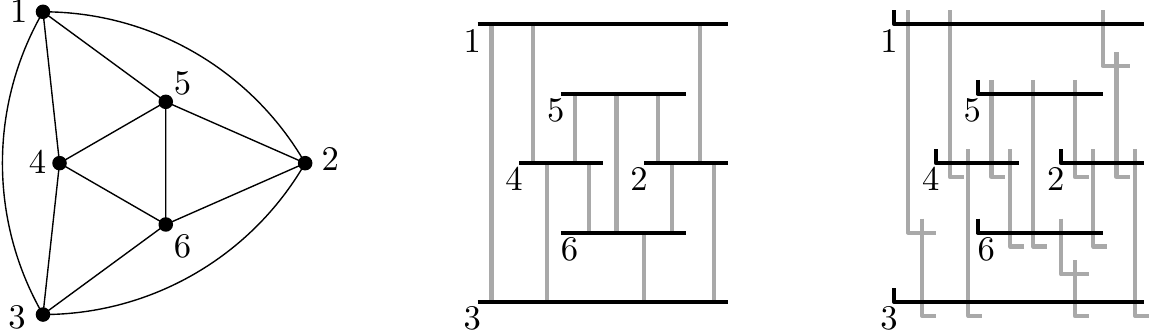}
    \caption{A planar graph $G$ on the left, a bar visibility representation of $G$ in the center, and an $\{\L\}$-representation of a full division of $G$ on the right. Here, the edges $\{1,2\}$, $\{1,3\}$ and $\{3,6\}$ are subdivided twice.}
    \label{fig:bar-visibility}
   \end{figure}
   
   It is now easy to interpret every segment as an \L, and replace a segment corresponding to an edge that is subdivided twice by two \L-shapes. Let us simply refer to Figure~\ref{fig:bar-visibility} again.
  
   \item Now assume that $G = (V,E)$ is any graph. We shall construct a $B_2$-VPG-representation of the complement $\bar{H}$ of $H = (V \cup W,E')$ with monotone increasing paths only. First, we represent the clique $\bar{H}[V]$. Let $V = \{v_1,\ldots,v_n\}$ and define for $i=1,\ldots,n$ the $2$-bend path $P(v_i)$ for vertex $v_i$ to start at $(i,0)$, have bends at $(i,i)$ and $(i+n,i)$, and end at $(i+n,n+1)$. See Figure~\ref{fig:B2-co-subdivision} for an illustration. For convenience, let us call these paths \emph{$v$-paths}.
   
   
   \begin{figure}[htb]
    \centering
    \includegraphics[scale=0.85]{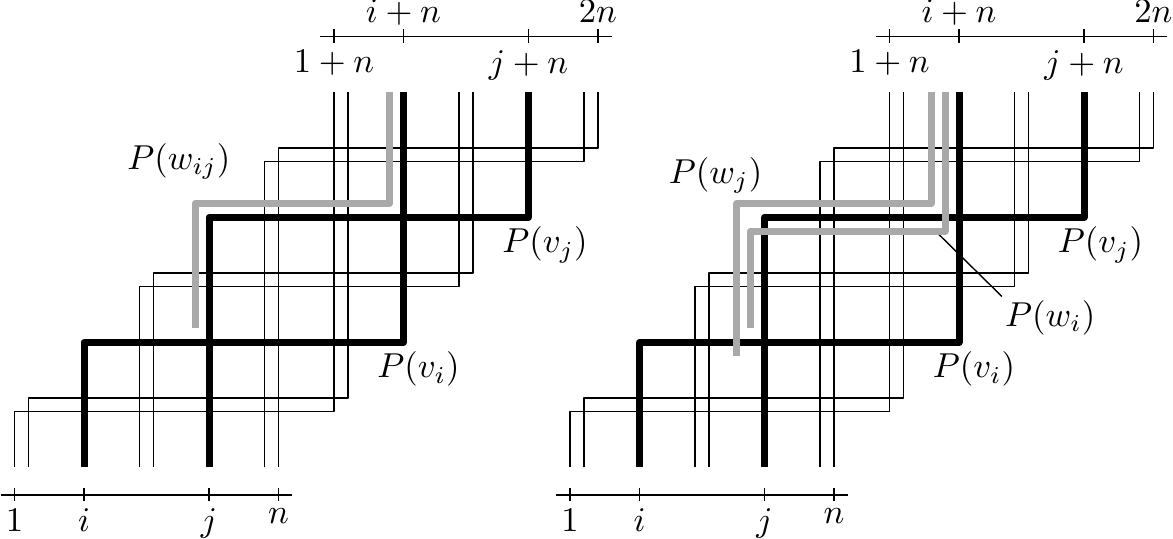}
    \caption{Left: Inserting the path $P(w_{ij})$ for a single vertex $w_{ij}$ subdividing the edge $v_iv_j$ in $G$. Right: Inserting the paths $P(w_i)$ and $P(w_j)$ for two vertices $w_i,w_j$ subdividing the edge $v_iv_j$ in $G$.}
    \label{fig:B2-co-subdivision}
   \end{figure}
   
   
   Next, we define for every edge of $G$ the $2$-bend paths for the one or two corresponding subdivision vertices in $\bar{H}$. We call these paths \emph{$w$-paths}. So let $\{v_i,v_j\}$ be any edge of $G$ with $i < j$. We distinguish two cases.
   
   
   \begin{enumerate}[leftmargin = 2.5em, label = \textit{Case }\arabic*.]
    \item The edge $\{v_i,v_j\}$ is subdivided by only one vertex $w_{ij}$ in $H$. We define the $w$-path $P(w_{ij})$ to start at $(j-\frac14,i+\frac14)$, have bends at $(j-\frac14,j+\frac14)$ and $(i+n-\frac14,j+\frac14)$, and end at $(i+n-\frac14,n+1)$, see the left of Figure~\ref{fig:B2-co-subdivision}.
    
    \item The edge $\{v_i,v_j\}$ is subdivided by two vertices $w_i,w_j$ with $\{v_i,w_i\}, \{v_j,w_j\} \in E(H)$. We define the start, bends and end of the $w$-path $P(w_i)$ to be $(j-\frac14,i+\frac14)$, $(j-\frac14,j-\frac14)$, $(i+n-\frac14,j-\frac14)$ and $(i+n-\frac14,n+1)$, respectively. The start, bends and end of the $w$-path $P(w_j)$ are $(j-\frac12,i-\frac14)$, $(j-\frac12,j+\frac14)$, $(i+n-\frac12,j+\frac14)$ and $(i+n-\frac12,n+1)$, respectively. See the right of Figure~\ref{fig:B2-co-subdivision}.
   \end{enumerate}
   
   
   \noindent
   It is easy to see that every $w$-path $P(w)$ intersects every $v$-path, except for the one or two $v$-paths corresponding to the neighbors of $w$ in $H$. Moreover, the two $w$-paths in Case~2 are disjoint. It remains to check that the $w$-paths for distinct edges of $G$ mutually intersect. To this end, note that every $w$-path for edge $v_iv_j$ starts near $(j,i)$, bends near $(j,j)$ and $(i+n,j)$ and ends near $(i+n,n)$. Consider two $w$-paths $P$ and $P'$ that start at $(j,i)$ and $(j',i')$, respectively, and bend near $(j,j)$ and $(j',j')$, respectively. If $j = j'$ then it is easy to check that $P$ and $P'$ intersect near $(j,j)$. Otherwise, let $j' > j$. Now if $j > i'$, then $P$ and $P'$ intersect near $(j',i)$, and if $j \leq i'$, then $P$ and $P'$ intersect near $(i+n,j')$.
   
   Hence we have found a $B_2$-\VPG-representation of $\bar{H}$, as desired. Let us remark, that in this representation some $w$-paths intersect non-trivially along some horizontal or vertical lines, i.e., share more than a finite set of points. However, this can be omitted by a slight and appropriate perturbation of endpoints and bends of $w$-paths.\qed
  \end{enumerate}
 \end{proof}

\section{Conclusions and Open Problems}
\label{sec:conclusions}

 Motivated by Middendorf and Pfeiffer's theorem (Theorem~\ref{thm:stretch-L-G} in~\cite{MP92}) that every~$\{\L,\G\}$-representation can be stretched into a \SEG-representation, we considered the question which subclasses of \SEG-graphs are actually $\{\L,\G\}$-graphs, or even $\{\L\}$-graphs. We proved that this is indeed the case for several graph classes related to planar graphs. We feel that the question whether $\PLANAR \subset \{\L,\G\}$, as already conjectured~\cite{CU13}, is of particular importance. After all, this, together with Theorem~\ref{thm:stretch-L-G}, would give a new proof for the fact that $\PLANAR \subset \SEG$.

%
%

 \begin{openproblem}
  Each of the following is open.
  \begin{enumerate}[label = (\roman*)]
   \item When can a $B_1$-\VPG-representation be stretched into a combinatorially equivalent \SEG-representation?
   \item Is $\{\L,\G\} = \SEG \cap B_1$-\VPG?
   \item Is every planar graph an $\{\L\}$-graph, or $B_1$-\VPG-graph?
   \item Does every planar graph admit an even subdivision whose complement is an $\{\L\}$-graph, or $B_1$-\VPG-graph?
   \item Recognizing $B_k$-\VPG-graphs is known to be NP-complete for each $k \geq 0$~\cite{CJKV12}. What is the complexity of recognizing $\{\L\}$-graphs, or $\{\L,\G\}$-graphs?
  \end{enumerate}
 \end{openproblem}


\noindent 
 Recall that $\cup_{k\geq 0} B_k$-$\VPG = \STRING$~\cite{ACGLLS12}. Chaplick \textit{et al.}~\cite{CJKV12} prove that $B_k$-\VPG{} $\subsetneq B_{k+1}$-\VPG{} for all $k\geq 0$ and also that $\SEG \nsubseteq B_k$-\VPG{} for each $k \geq 0$, even if \SEG{} is restricted to three slopes only. Another natural subclass of \STRING, which is in no inclusion-relation with \SEG, is the class \COCO{} of co-comparability graphs~\cite{golumbic1983coco}. However, one can prove a result similar to the previous one concerning $B_k$-\VPG-graphs and \STRING-graphs:

 \smallskip

There is no $k \in \mathbb{N}$ such that $B_k$-$\VPG \supset
\COCO$, compare also Figure~\ref{fig:graph-classes}. A proof can be given along the ``degrees of freedoms''
approach of Alon and Scheinerman~\cite{AS88}, i.e., by
counting the graphs in the respective sets.

First, Alon and Scheinerman consider the class $P(n,t)$
of $t$-dimensional posets on $n$ elements and show that
for fixed $t$ the growth of $\log |P(n,t)|$ behaves like $nt\log n$.
If $CC(n,t)$ denotes the class of cocomparability graphs of posets in $P(n,t)$, then with easy adaptations we obtain $\log |CC(n,t)| \geq n(t-1-o(1))\log n.$

On the other hand, every path of a $B_k$-\VPG-representation can be encoded by $k+4$ numbers.
Whether two paths intersect can be answered by looking at the signs of few low degree polynomials in $2k+8$ variables evaluated at the encodings of the two paths, meaning that the class $B_k$-\VPG{} has $k+4$ degrees of freedom.
Alon and Scheinerman show how to use Warren's Theorem~\cite{warren1968} to get an upper bound on the size of such a class.
Indeed, the logarithm of the number of $B_k$-\VPG{}-graphs on $n$ vertices is at most $O(1)nk\log n$. 
 
Comparing the numbers we find that there is a cocomparability graph of a $(k+2)$-dimensional poset that is not a $B_k$-$\VPG$-graph.
On the other hand it is easy to find a $B_k$-$\VPG$-representation for the cocomparability graph of any given $(k+1)$-dimensional poset.
So we have for every $k \in \mathbb{N}$
\[
 \bigcup_{n \in \mathbb{N}}CC(n,k+1) \subseteq B_k\text{-}\VPG{} \nsupseteq \bigcup_{n \in \mathbb{N}}CC(n,k+2).
\]
Recently, this has been proved independently by Cohen, Golumbic, Trotter and Wang~\cite{Coh-14} and hence we omit more details here.

Comparing this to the result that there is no $k$ such that $\SEG \subseteq B_k$-\VPG{} a natural question arises. Is there a parameter of \SEG-graphs or \SEG-representations ensuring few bends in their \VPG{}-representations?
It is known that the number of slopes in the \SEG-representation is not the right answer~\cite{CJKV12}.
 
 
\bibliographystyle{abbrv}
\bibliography{lit}

\end{document}